\documentclass[]{birkmult}
\usepackage{soul,cancel}

 \usepackage{color}

\usepackage{amsmath}
\usepackage{amsfonts}
\usepackage{amssymb}
\usepackage{color}
\usepackage{srcltx}

 \newtheorem{thm}{Theorem}
 \newtheorem{cor}[thm]{Corollary}
 \newtheorem{lem}[thm]{Lemma}

 \newtheorem{rem}[thm]{Remark}
 
\newcommand{\eq} [1] {\begin{equation}\label{#1}\quad}
\newcommand{\en} {\end{equation}}
\newcommand{\scal}[1]{\langle#1\rangle}
\newcommand{\norm}[1]{\left\Vert#1\right\Vert}

\newcommand{\Int}{\operatorname{int}}
\newcommand{\tr}{\operatorname{tr}}
\renewcommand{\Re}{\operatorname{Re}}

\newcommand{\C}{\mathbb{C}}
\newcommand{\R}{\mathbb{R}}

\def\cC{\mathcal C}
\def\pt{\partial}
\def\CR{\color{red}}

\begin{document}
%
%

\title{The possible  shapes of numerical ranges}

\author[Helton]{J. William Helton}
\address{Department of Mathematics\\
University of California San Diego\\
9500 Gilman Drive \\ La Jolla, CA 92093-0112\\
USA} \email{helton@math.ucsd.edu} \thanks{Research
of the first author supported by NSF grants DMS-0700758,
DMS-0757212, and the Ford Motor Co.}

\author[Spitkovsky]{I.M. Spitkovsky}
\address{Department of Mathematics\\ College of William and Mary\\
Williamsburg, VA 23187\\ USA} \email{ilya@math.wm.edu,
imspitkovsky@gmail.com}



\begin{abstract}
Which convex subsets of $\C$ are the numerical range $W(A)$ of
some matrix $A$? This paper gives a precise characterization of these
sets. In addition to this we show that  for any $A$ there
exists a symmetric $B$ of the same size such that $W(A)=W(B)$
thereby settling an open question
from \cite{CLT11}.
\end{abstract}

\maketitle

Consider $\C^d$,  the standard complex inner product space. Let
$\scal{\cdot,\cdot}$ denote its scalar product, and $\norm{\cdot}$
the related norm. The {\em numerical range}  $W(A)$ of a $d\times d$
matrix $A$ is defined as \eq{nr} W(A)=\{ \scal{Ax,x}\colon
\norm{x}=1\}.  \en  It is well known that $W(A)$ is a compact convex
subset of $\C$ containing the spectrum of $A$; see, e.g., monographs
\cite{GusRa,HJ2} for these and other properties, as well as for the
history of the subject. In this short note  we give an answer to the
question of exactly which sets $W$ actually are the numerical range
of some matrix $A$.

 This  question was originally raised in
 Kippenhahn's 1951 article
 \cite{Ki} (see also a more accessible English translation \cite{Ki08})
which gave several non-trivial necessary conditions on
 the ``geometrical shape" of a numerical range.

However, a necessary and sufficient condition remained
open\footnote{We would like to thank P. Y. Wu for
discussion of this issue during XXVII South Eastern Analysis Meeting
in Gainesville, FL.}.
 One can be obtained by
the observation that  curves critical to the problem were effectively
classified in \cite{HeVin07}. Didier Henrion in \cite{Hen10}  makes
such a connection\footnote{We are especially grateful to Bernd
Sturmfels for bringing \cite{Hen10} to our attention.} and more, and
states explicitly one side (necessary) of the characterization of
numerical range. While all components of our paper  can easily be
extracted from \cite{Hen10} by hose comfortable with the theory in
\cite{HeVin07} , we think our short note will nevertheless be useful to
the numerical range community, at least for expository purposes.
 In particular, our
Theorem \ref{main} explicitly states a necessary and sufficient
condition.


Our characterization of numerical ranges is in terms of a type of dual
convex set. { For any  set $S\subset\R^n$} its {\em polar}
is defined as \eq{pol} S_*=\{ x\in\R^n\colon \sup_{y\in S}\,
\scal{x,y}\leq 1\} \en
 (see, e.g., \cite{BeNe01,Rock97}).
 The set $S_*$ is closed, convex,  and contains 0.
 Clearly (see also \cite[Corollary 14.5.1]{Rock97}),
 0 is an interior point of $S_*$ if and only if $S$ is bounded. If $S$ itself is closed, convex and contains 0, then
 \eq{dual} (S_*)_*=S
\en  \cite[Theorem  14.5]{Rock97}.

The next result
provides  an explicit description of polar sets of numerical ranges.
In some form it goes back many years, at least to \S 3 \cite{Ki}.
%
 A different point of view (in a more
general setting) is presented in \cite[Section 5]{RoSt11} (there the
term dual is used in place of polar).

\begin{lem}
\label{l:dual}
Let $A\in\C^{d\times d}$.
Then
\eq{lmirep}
 W(A)_*= \{z=\xi + i\eta \colon I - \xi H-\eta K \text{ is positive semi-definite} \}. \en
Here $H$ and $K$ are hermitian matrices from the representation
\eq{AHK} A=H+iK.\en
 \end{lem}

\begin{proof}
Directly from the definitions (\ref{nr}) and (\ref{pol}) it follows that
\begin{multline*}
 W(A)_*=\{ z\colon
\Re (\scal{Av,v}\overline{z})\leq 1 \text{ for all } v\in\C^d\text{ with }
\norm{v}=1\} \\ = \{ z\colon \scal{(\Re(\overline{z}A)v,v})\leq 1 \text{
for all } v\in\C^d\text{ with } \norm{v}=1\} \\ = \{ z\colon
I-\Re(\overline{z}A)\text{ is positive semi-definite} \}\\ = \{
\xi+i\eta\colon I-\xi H-\eta K\text{ is positive semi-definite}\}.
\end{multline*}
\end{proof}

Common terminology is that \eqref{lmirep} is a {\em linear matrix
inequality  {\em (LMI for short)} representation} for $W(A)_*$ and the
lemma says that if a set $W\subset \C$ is a numerical range, then its
polar has an LMI representation. The paper \cite{HeVin07} describes
precisely the sets $\cC$ in $\R^2$, hence in $\C$, which have an LMI
representation. It characterizes them as ``rigidly convex" a term we set
about to define. An {\it algebraic interior} $\cC$  has a defining
polynomial $q$, namely $\cC$ is the closure of the connected
 component of $\cC:=\{z: \ q(z)>0\}$
containing 0. A minimum degree defining polynomial for $\cC$ is
unique (up to a constant), see Lemma 2.1 \cite{HeVin07}
 and its degree we call the degree of $\cC$.
A convex set $\cC$ is called {\it rigidly convex} provided it is an
algebraic interior  and it has a defining polynomial
$q$ which satisfies the {\it real zero} (RZ) condition, namely,
$$if \ \mu \in \C\ and \ q(\mu z)=0, \ then \ \mu\in\R.$$
Our main theorem is:
\begin{thm}\label{main}
A subset $W$ of $\C$ is the numerical range of some $d\times d$
matrix $A$ if and only if its polar $W_* $ is rigidly convex of degree less
than or equal to $d$.
\end{thm}

\begin{proof}
Given $A=H+i K$, observe
that $p$ defined by
\begin{equation}
\label{eq:det}
p(z)=\det(I-\xi H-\eta K)
\end{equation}
 is an $RZ$  polynomial,
since all eigenvalues of a symmetric matrix are real.
Moreover, $W(A)_*$ coincides with the
closure of the connected component of $\{z\colon p(z)>0\}$
containing zero. Thus the set $W(A)_*$ is
rigidly convex.

However, Theorem 3.1 of \cite{HeVin07} says that converse also
holds \footnote{For perspective, \cite{LPR05} showed that  the proof of Theorem 3.1 in \cite{HeVin07}
implies a 1958 conjecture of Peter Lax is true. In this context we might
describe the characterization of numerical ranges (Theorem \ref{main}) as ``polar" to the Lax Conjecture.}
:
if $V$ is rigidly convex, then there exist real
symmetric matrices $H,K$ such that
 \eq{Ws}  V=\{ z=\xi+i\eta\colon
I-\xi H-\eta K \text{ is positive semi-definite}\}.
\en
 Consequently,
$V=W(B)_*$ for $B=H+iK$.
Moreover, we can do this
with an  $H,K$ whose dimension is the degree of $V$.
\end{proof}
The forward side of Theorem \ref{main} is in \cite{Hen10}
 (stated in the language of homogeneous coordinates,
and emphasizing that numerical ranges are affine projections of
semi-definite cones). The converse follows easily from ingredients
there,  though it is not stated  explicitly.

Note that the matrix $B$
constructed in the proof of Theorem~\ref{main}
is symmetric along with $H,K$.
This yields an affirmative answer to the
 question stated in \cite{CLT11} (raised by the referee of
the latter):

\begin{cor} For every
$d\times d$ matrix $A$ there exists a symmetric $d\times d$ matrix
$B$ such that $W(B)=W(A)$.
\end{cor}

Duality (\ref{dual}) allows us to restate Theorem~\ref{main}
in the following form.

\begin{cor}
A subset $W$ of $\C$ is the numerical range of some $d\times d$
matrix $A$ if and only if it is a translation of the polar of a rigidly
convex set of degree less than or equal to $d$.
\end{cor}

\begin{proof}
For a given  $d\times d$ matrix
$A$, pick $\lambda
\in W(A)$ and let $A_0=A-\lambda I$.
 By
Theorem~\ref{main}, the polar set $V$ of $W(A_0)$  is rigidly convex
and has  degree not exceeding $d$.
But $0\in W(A_0)$, so that due to
(\ref{dual}) we have $W(A_0)=V_*$.
Consequently, $W(A)=W(A_0)+\lambda$ is a
translation of $V_*$

Conversely, if $W$ is a translation of $V_*$ for some rigidly convex set
$V$ of degree not exceeding $d$, then $W-\lambda=V_*$ for some
$\lambda\in\C$. Applying (\ref{dual}) to $S=V$, we conclude that
$(W-\lambda)_*=V$. By Theorem~\ref{main}, $W-\lambda=W(A_0)$ for
some $d\times d$ matrix $A_0$, so that $W=W(A_0+\lambda I)$.
\end{proof}

\begin{rem}\rm
If the matrices $H,K$ from representation (\ref{AHK}) are linearly dependent with $I$,
then the set $V$ in (\ref{Ws}) is unbounded.
Moreover, $V$ stays unbounded
under translations of $A$. In other words, $W(A)$ in this case has
empty interior. This agrees with the fact that $A$ in this (and only this)
case has the form $\alpha R+\beta I$ for some hermitian $R$ and
$\alpha,\beta\in\C$, and $W(A)$ is therefore a (closed) line segment.
In all other cases the interior of $W(A)$ is non-empty, and
$W(A-\lambda I)_*$ is bounded for any
$\lambda$ lying in the interior of $W(A)$.
One such value of $\lambda$ is  $\lambda=\tr (A)/d$.
\end{rem}

\providecommand{\bysame}{\leavevmode\hbox
to3em{\hrulefill}\thinspace}
\providecommand{\MR}{\relax\ifhmode\unskip\space\fi MR }
\providecommand{\MRhref}[2]{%
  \href{http://www.ams.org/mathscinet-getitem?mr=#1}{#2}
} \providecommand{\href}[2]{#2}

\end{document}